\documentclass[12pt]{amsart}
\usepackage{amsfonts}
\usepackage{amssymb}
\usepackage{amscd}

\newtheorem{claim}{Claim}
\newtheorem{cor}{Corollary}
\newtheorem{duffy}{Definition}
\newtheorem{lemma}{Lemma}
\newtheorem{prop}{Proposition}
\newtheorem{thm}{Theorem}

\begin{document}

\title{Linear systems attached to cyclic inertia}
\author{Marco A. Garuti}
\address{\hskip-\parindent
        Marco A. Garuti\\
        Dipartimento di Matematica Pura ed Applicata\\
        Universit\`a degli Studi di Padova\\
        Via Belzoni 7, 35131\\
        Padova - ITALY}
\email{mgaruti@galileo.math.unipd.it}

\thanks{Research at MSRI is supported in part by 
NSF grant DMS-9701755.}

\begin{abstract}
We construct inductively an equivariant compactification
of the algebraic group
${\mathbb W}_n$ of Witt vectors of finite length over a
field of characteristic $p>0$.
We obtain smooth projective rational varieties
$\overline{\mathbb W}_n$, defined over $\mathbf F_p$;
the boundary is a divisor whose reduced subscheme has
normal crossings.

The Artin-Schreier-Witt isogeny $F-1:{\mathbb W}_n\to
{\mathbb W}_n$ extends to a finite cyclic cover 
${\mathbf\Psi}_n:\overline{\mathbb W}_n\to
\overline{\mathbb W}_n$ of degree $p^n$ ramified at the
boundary.
This is used to give an extrinsic description of the
local behavior of a separable cover of curves 
in char. $p$ at a wildly ramified point whose inertia
group is cyclic.

In an appendix, we give an elementary computation of the
conductor of such a covering, which can otherwise be
determined using class field theory.
\end{abstract}

\maketitle


Let $f:C\to D$ be a finite separable cover 
of (germs of) curves 
over a perfect field $k$ of 
characteristic $p>0$.
Suppose that $f$ is 
Galois of group $G$ and let 
$y_0\in C$ be a ramification point
at which the residue field extension 
$k(y_0)/k(x_0)$ is separable, 
where $x_0=f(y_0)$. 
Recall (\cite{CL}, chap. IV), 
that the higher ramification
groups 
$G_i$ (in lower numbering) are defined as 
$$
G_i=\left\{g\in G|\, 
v_{y_0}\left(g(t)-t\right)\geq i+1\right\}
$$
where $t$ is a local parameter at $y_0$; 
they do not depend on the choice of $t$.
The inertia group $G_0$ is an extension of 
a $p$-group by a cyclic group of order 
prime to $p$ and, for $i>0$, 
the $G_i$ are $p$-groups;
they are trivial for 
large values of $i$.
The degree (at $y_0$) of the different
$\mathfrak{D}_{C/D}=Ann\,(\Omega_{C/D}^1)$ 
of the extension $C/D$ can be expressed as
$$
\deg\mathfrak{D}_{C/D,y_0}
=\sum_{i=0}^\infty\left(|G_i|-1\right)
$$

Using Herbrand's $\varphi$ and $\psi$ functions, 
one can introduce (\emph{loc.cit.}, \S 3) higher
ramification groups in upper numbering
$G^{\varphi(r)}=G_r$, indexed over 
$\mathbf{R}_{\geq0}$.
Hence, if $e$ is the largest integer such that 
$G_e$ is nontrivial, the real number 
$$
m=\varphi(e)
$$
is the largest for which $G^m$ is nontrivial. 
When $G$ is abelian,
the classical Hasse-Arf theorem (\emph{loc.cit.}) 
guarantees that this number is in fact an integer. 
We shall write $m(f)$
and $e(f)$ to emphasize that these integers are 
related to the cover $f:C\to D$. 
The conductor (at $y_0$) of 
the extension $C/D$ is defined as
$$
\mathrm{cond}\,(C/D,y_0)=m(f)+1.
$$

Geometrically, the conductor can be characterized 
as follows 
(\cite{gacc}, chap. VI, 2.12 and \cite{clac}, 3.6): 
$m+1=\deg_{y_0}\mathfrak{m}$, where $\mathfrak{m}$ 
is the smallest modulus
such that the cover $C/D$ is isomorphic to a 
pull-back from an isogeny of
the generalized jacobian $J_{\mathfrak{m}}$ of $D$.

\bigskip

We shall be concerned with the situation 
where $G_0$ is cyclic of order $p^n$.
By Artin-Schreier-Witt theory, there are 
\'etale neighborhoods $V$ of $y_0$ and 
$U$ of $x_0$ and a map $\underline u$, 
from $U^\prime=U-\{x_0\}$ to the group 
${\mathbb W}_n$ of Witt vectors of length $n$,
such that $V^\prime=V-\{y_0\}$ is the 
fibred product
$$
\CD
V^\prime @>>> {\mathbb W}_n \\
@V{f|_{V^\prime}}VV    @VV{F-1}V \\
U^\prime @>{\underline u}>> {\mathbb W}_n
\endCD
$$
The maps $\underline  u$ and
$\underline u+F\underline w-\underline w$, 
define isomorphic covers
for any Witt-vector valued function
$\underline w$.

\bigskip

We will make the following 
assumptions: $f$ itself is cyclic of degree $p^n$, 
totally ramified at $y_0$; furthermore, 
$U=D$ and $V=C$ and $\underline u$ 
is defined over $k$. 
Working locally for the \'etale topology, 
this is no great loss of generality.

\bigskip
We are going to show that, 
for every positive integer $n$,
there is a smooth projective rational 
variety ${\overline{\mathbb W}}_n$
over the prime field $\mathbf{F}_p$, 
equipped with a tautological line bundle 
$\mathcal{O}_{{\overline{\mathbb W}}_n}(1)$
such that:
\begin{enumerate}
\item ${\overline{\mathbb W}}_n$ contains 
${{\mathbb W}}_n$ as an affine open dense subvariety; 
the boundary $B_n={\overline{\mathbb W}}_n
\setminus{{\mathbb W}}_n$ is a divisor such that
$\mathcal{L}(B_n)=
\mathcal{O}_{{\overline{\mathbb
W}}_n}(1)$  and $(B_n)_{\mathrm{red}}$ has normal
crossings.
\item The isogeny $F-1$ extends to a separable 
cyclic cover $\mathbf\Psi_n:
\overline{\mathbb W}_n\to\overline{\mathbb W}_n$ 
branched along $B_n$.
\end{enumerate}

\noindent
Our main result is then:

\begin{thm}
Let $f:C_n\to D$ be a separable cyclic 
cover of degree $p^n$ of curves
over a perfect field $k$, totally 
branched above $x_0\in D(k)$.
Let $\underline u=
(u_0,\dots,u_{n-1})$ 
be a Witt vector of rational 
functions on $D$, with poles 
of order $\nu_i$, prime to $p$,
 at $x_0$, such that 
$k(C_n)/k(D)$ is defined by 
$F(\underline{Y})-\underline Y
=\underline u$.

Then $\underline u$ extends to a morphism
$\underline u: D\to \overline{\mathbb W}_n\times k$ such 
that $C_n$ is isomorphic to the normalization 
of the fibred product 
$D\times_{\overline{\mathbb W}_n}\overline{\mathbb W}_n$
and
$$
\mathrm{cond}\,(C_n/D,y_0)-1
=\deg\underline u^\ast
\mathcal{O}_{{\overline{\mathbb W}}_n}(1)=
\max_{0\leq i\leq n-1}\{p^{n-i-1}\nu_i\}. 
$$
\end{thm}

\bigskip

This work was completed while visiting the M.S.R.I.
(october-december 1999) during the program on
{\it Galois Groups and Fundamental Groups.}
I would like to thank the participants for many
stimulating conversations and in particular Ahmed Abbes
and Mohamed Sa\"\i di for pointing out to me the
references \cite{bry} and \cite{schmid} respectively.

\bigskip
\bigskip

We begin with the following
\begin{duffy}
For every positive integer $n$, 
define inductively a 
variety ${\overline{\mathbb W}}_n$
over the prime field $\mathbf{F}_p$, 
equipped with a tautological line bundle 
$\mathcal{O}_{{\overline{\mathbb W}}_n}(1)$, 
by letting
$$\begin{array}{rl}
\left({\overline{\mathbb W}}_1,
\mathcal{O}_{{\overline{\mathbb W}}_1}(1)\right)
& = \left(\mathbf{P}^1,
\mathcal{O}_{\mathbf{P}^1}(1)\right)\\
\left({\overline{\mathbb W}}_{n+1},
\mathcal{O}_{{\overline{\mathbb W}}_{n+1}}(1)\right)
& = \left(\mathbf{P}\left(
\mathcal{O}_{{\overline{\mathbb W}}_n}
\oplus\mathcal{O}_{{\overline{\mathbb W}}_n}(p)\right),
\mathcal{O}_{\mathbf P}(1)\right)
\end{array}
$$
(as usual, $\mathbf P(\mathcal E)$ denotes 
the projective bundle
associated to a locally free module 
$\mathcal E$ and 
$\mathcal{O}_{\mathbf P}(1)$ its 
tautological line bundle 
cf.\cite{EGA2}, 4.1.1; for any $m\in\mathbf Z$, 
put $\mathcal{O}(m)=\mathcal{O}(1)^{\otimes m}$).
\end{duffy}

By definition, there is a natural projection
\begin{equation}
\label{res}
\mathbf{r}:
{\overline{\mathbb W}}_{n+1}\longrightarrow{\overline{\mathbb W}}_n
\end{equation}
which is a fibration by projective lines;
it follows immediately that 
${\overline{\mathbb W}}_n$ is a smooth projective rational 
variety of dimension $n$.
Note that $\mathcal{O}_{{\overline{\mathbb W}}_{n}}(1)$ is
not ample (for $n\geq2$); 
it is only ample relative to $\mathbf r$.

On the other hand, ${\overline{\mathbb W}}_{n+1}$ 
is the projective closure of the vector bundle
$\mathbf V(\mathcal{O}_{{\overline{\mathbb W}}_{n}}(p))$
over ${\overline{\mathbb W}}_n$ (cf.~\cite{EGA2}, 8.4) and 
there are two canonical sections of $\mathbf r$. 
The zero section, whose divisor we will 
denote $Z_{n+1}$, arises from the 
exact sequence
$$
0\longrightarrow 
\mathcal{O}_{{\overline{\mathbb W}}_n}
\longrightarrow 
\mathcal{O}_{{\overline{\mathbb W}}_n}
\oplus\mathcal{O}_{{\overline{\mathbb W}}_n}(p)
\longrightarrow 
\mathcal{O}_{{\overline{\mathbb W}}_n}(p)
\longrightarrow 0
$$
which gives a global section of
$\mathcal{O}_{{\overline{\mathbb W}}_n}\oplus
\mathcal{O}_{{\overline{\mathbb W}}_n}(p)=
\mathbf r_\ast
\mathcal{O}_{{\overline{\mathbb W}}_{n+1}}(1)$,
whence a global section of
$\mathcal{O}_{{\overline{\mathbb W}}_{n+1}}(1)$
(the latter restricts to $\mathcal{O}_{\mathbf P^1}(1)$
on the fibres of $\mathbf r$, 
so its higher direct images vanish and 
$H^i(\overline{\mathbb W}_n,\mathbf r_\ast
\mathcal{O}_{{\overline{\mathbb W}}_{n+1}}(1))
=H^i(\overline{\mathbb W}_{n+1},
\mathcal{O}_{{\overline{\mathbb W}}_{n+1}}(1))$
for all $i$).

The infinity section is obtained by the same construction 
composed with the isomorphism 
$$
\mathbf{P}\left(\mathcal{O}_{{\overline{\mathbb W}}_n}
\oplus\mathcal{O}_{{\overline{\mathbb W}}_n}(p)\right)
\simeq \mathbf{P}\left(\mathcal{O}_{{\overline{\mathbb W}}_n}
(-p)\oplus\mathcal{O}_{{\overline{\mathbb W}}_n}\right).
$$
We denote $\Sigma_{n+1}$ the corresponding divisor
in $H^0(\overline{\mathbb W}_{n+1},
\mathcal{O}_{{\overline{\mathbb W}}_{n+1}}(1)\otimes
{\mathbf r}^\ast\mathcal{O}_{{\overline{\mathbb W}}_n}(-p))$.

\bigskip

As suggested by the notation, ${\overline{\mathbb W}}_n$ is a 
compactification of the Witt group scheme ${{\mathbb W}}_n$; 
more precisely:

\begin{prop}
There is a natural system of open immersions
$j_n:{{\mathbb W}}_n\subset{\overline{\mathbb W}}_n$, compatible w.r.t. 
the projections $\mathbf r$ and such that 
$\mathcal{O}_{{\overline{\mathbb W}}_n}(1)$ 
trivializes over $j_n({{\mathbb W}}_n)$.

The closed subscheme 
$B_n={\overline{\mathbb W}}_n\setminus j_n({{\mathbb W}}_n)$ 
is a divisor linearly equivalent to
$\Sigma_n+p\ \mathbf r^\ast B_{n-1}$ 
(and hence to $Z_n$)
and $(B_n)_{\mathrm{red}}$ 
has normal crossings.
\end{prop}
Proof: Indeed, we have a canonical embedding 
$j_1:\Bbb G_{\mathrm a}\cong\mathbf A^1\subset\mathbf P^1$; 
by induction, assume $\mathcal{O}_{{\overline{\mathbb W}}_n}(1)$ 
trivializes over the affine open subset $j_n({\mathbb W}_n)=U_n$: 
then, by definition (cf.\cite{EGA2}, 3.1), 
$\mathbf{r}^{-1}(U_n)=\mathbf P_{U_n}^1$.
Clearly, $\mathcal{O}_{{\overline{\mathbb W}}_{n+1}}(1)$
trivializes over the complement $U_{n+1}$ 
of the infinity section, 
and, identifying ${\mathbf r}^{-1}(0)\cap U_{n+1}$ with 
$\mathbf A^1$, we can choose an isomorphism 
$j_{n+1}:{\mathbb W}_{n+1}\cong
U_{n+1}$  fitting in the diagram
$$
\CD
\mathbf A^1
@ >>> U_{n+1} @>{{\mathbf r}|_{U_{n+1}}}>> U_n \\ 
@Aj_1AA   @AAA    @AAj_n A \\
\Bbb G_{\mathrm a} @>{V^n}>>{\mathbb W}_{n+1} 
@>{\mathbf r}>> {\mathbb W}_n
\endCD
$$
where the bottom one is the exact sequence of group 
schemes defining ${\mathbb W}_{n+1}$.

\noindent
The second statement follows easily, as 
$j_{n+1}({\mathbb W}_{n+1})$ has been defined 
as the intersection of $\mathbf{r}^{-1}(U_n)$, 
which is the complement of the vertical divisor 
$\mathbf r^\ast(pB_n)$, with the complement of 
the horizontal infinity divisor $\Sigma_{n+1}$.

\begin{cor}
\label{bord}
For all $1\leq i\leq n$, let 
$B_{n,i}=({\mathbf r}^{n-i})^\ast(\Sigma_i)$
be the pull back of the infinity divisor 
on $\overline{\mathbb W}_{i}$
via the iteration of the 
projection (\ref{res}). 
The boundary divisor is then
$$
B_n=\sum_{i=1}^np^{n-i}
B_{n,i}
$$
\end{cor}

From now on, we identify ${{\mathbb W}}_n$ with its 
image in ${\overline{\mathbb W}}_n$.
Choose a system of variables 
$Y_0,Y_1,\dots,Y_n,\dots$ such that 
${{\mathbb W}}_n=\mathrm{Spec}\ 
\mathbf F_p[Y_0,\dots,Y_{n-1}]$ and 
the projection ${\mathbf r}$ 
corresponds to the inclusion
$\mathbf F_p[Y_0,\dots,Y_{n-1}]\subset
\mathbf F_p[Y_0,\dots,Y_{n}]$.
On the ring of polynomials
$$
{\mathsf H}=\mathbf F_p[T,Y_0,\dots,Y_n,\dots]
$$
we define a grading by $\deg T=1$ and $\deg Y_i=p^i$
(recall \cite{AC} that with this choice for the degrees, 
the polynomials defining the sum, product etc. 
of Witt vectors are isobaric).
As usual, we denote by ${\mathsf H}_d$ the $d$-th graded
piece, consisting of homogeneous elements of degree $d$.

\begin{lemma}
\label{H=Gamma}
The tautological line bundle
$\mathcal{O}_{{\overline{\mathbb W}}_n}(1)$ 
is generated by global sections;
there is a canonical isomorphism
$H^0(\overline{\mathbb W}_n,
\mathcal{O}_{{\overline{\mathbb W}}_n}(1))
={\mathsf H}_{p^{n-1}}$ 
under which $Z_n=(Y_{n-1})_0$
and $B_n=(T^{p^{n-1}})_0$.
\end{lemma}
Proof: The assertions being clear for $n=1$, 
the proof is by induction.
We need a classical explicit description of the 
projective bundle $\overline{\mathbb W}_{n+1}$.

\noindent
By inductive assumption, the complete linear 
series $|pZ_n|$ defines a 
morphism
$$
\phi_{|pZ_n|}:
\overline{\mathbb W}_{n}\longrightarrow
\mathbf P^m=\mathbf PH^0(\overline{\mathbb W}_n,
\mathcal{O}_{{\overline{\mathbb W}}_n}(p)).
$$
The vector space decomposition 
${\mathsf H}_{p^n}=H^0(\overline{\mathbb W}_n,
\mathcal{O}_{{\overline{\mathbb W}}_n}(p))
\oplus\langle Y_n\rangle$
induces a rational map 
$$
\pi:\mathbf P^{m+1}=
\mathbf P{\mathsf H}_{p^n}
\dasharrow\mathbf P^m
$$
defined outside $[0:\dots:0:1]$. The indeterminacy 
is solved on the blow-up of this point, 
$\tilde\mathbf P^{m+1}=
\mathbf P(\mathcal O_{\mathbf P^m}\oplus
\mathcal O_{\mathbf P^m}(1))$, and 
$\overline{\mathbb W}_{n+1}$ is canonically 
isomorphic to the fibred product: 
$$
\CD
\overline{\mathbb W}_{n+1} @>>> \tilde\mathbf P^{m+1} \\
@V{\mathbf r}VV       @VV{\tilde\pi}V \\
\overline{\mathbb W}_{n} @>\phi_{|pZ_n|}>> \mathbf P^m
\endCD
$$
To complete the proof of the lemma it suffices to 
remark that the composite morphism 
$\overline{\mathbb W}_{n+1}\to\tilde\mathbf P^{m+1}
\to\mathbf P^{m+1}$
is given by the complete linear system
$|Z_{n+1}|$; 
its image is the projective cone
over 
$\phi_{|pZ_n|}(\overline{\mathbb W}_{n})$.

\bigskip

\begin{prop}
The action of ${{\mathbb W}}_n$ on itself extends to 
an action of ${{\mathbb W}}_n$ on $\overline{\mathbb W}_n$
 and the line bundle 
$\mathcal{O}_{{\overline{\mathbb W}}_n}(1)$
is stable under this action.

In particular, 
$\mathbf Z/p^n\mathbf Z={{\mathbb W}}_n(\mathbf F_p)$ 
acts on $\overline{\mathbb W}_n$ 
and the 
the Artin-Schreier-Witt isogeny 
$(F-1):{\mathbb W}_{n} \to {{\mathbb W}}_n$
extends to a cyclic cover 
$$
\mathbf\Psi_n: {\overline{\mathbb W}}_{n}
\longrightarrow{\overline{\mathbb W}}_{n}
$$
of degree $p^n$, defined over $\mathbf F_p$,
commuting with the projection $\mathbf r$. 
It is branched along the boundary divisor 
$B_n$ and the inertia subgroup of the infinity 
section $\Sigma_{n}$ has order $p$.

The induced map $\mathbf\Psi_n^\ast$ on 
the group of $k$-cycles modulo
rational equivalence is multiplication by $p^{n-k}$. 
\end{prop}
Proof: Clearly, $\Bbb G_{\mathrm a}$ acts on 
$\mathbf P^1$ and the isogeny
$\wp:\Bbb G_{\mathrm a}\to\Bbb G_{\mathrm a}$ 
extends to a cyclic cover of degree $p$ of 
$\mathbf P^1$ 
branched at $\infty$; we proceed by induction.

\noindent
The group ${\mathbb W}_{n+1}$ acts on itself 
by translations;
recall (e.g. \cite{AC})
that the sum of two Witt vectors 
$\underline{X},\underline{Y}$ is
expressed by polynomials 
$S_i(\underline{X};\underline{Y})\in\mathbf{Z}
\left[X_0,\dots,X_i,Y_0,\dots,Y_i\right]$ 
characterized by the relations
$$
\Phi_n(S_0,\dots,S_n)=
\Phi_n(X_0,\dots,X_n)+\Phi_n(Y_0,\dots,X_n)
$$
If the variables $X_i$ and $Y_i$ are given 
the weight $p^i$, then $S_i$
is isobaric of weight $p^i$. Moreover:
\begin{equation}
\label{S_n}
S_i(\underline{X};\underline{Y})=
X_i+Y_i+c_i(\underline{X};\underline{Y}),
\qquad c_i\in\mathbf{Z}\left[X_0,
\dots,X_{i-1};Y_0,\dots,Y_{i-1}\right]
\end{equation}
Therefore, given an $\mathbf F_p$-algebra $A$,
for any 
$\underline a\in{\mathbb W}_{n+1}(A)$, the polynomial
$$
T^{p^n}S_n(\frac{Y_0}{T},\dots,\frac{Y_n}{T^{p^n}};
\underline a)=Y_n+a_nT^{p^n}
+T^{p^n}c_n(\frac{Y_0}{T},
\dots,\frac{Y_{n-1}}{T^{p^{n-1}}};\underline a)
$$
is a homogeneous element in 
${\mathsf H}_{p^n}\otimes A$.

\noindent
As  
$\mathcal{O}_{{\overline{\mathbb W}}_n}(1)$ is stable under 
the action, ${\mathbb W}_{n}(A)$
acts $A$-linearly on
$H^0(\overline{\mathbb W}_{n,A},
\mathcal{O}_{{\overline{\mathbb W}}_n}(p))$: 
we extend it to an action of 
${\mathbb W}_{n+1}(A)$ on 
${\mathsf H}_{p^n}\otimes A$ by 
\begin{equation}
\label{act}
\underline a\cdot Y_n=Y_n+a_nT^{p^n}
+T^{p^n}c_n(\frac{Y_0}{T},
\dots,\frac{Y_{n-1}}{T^{p^{n-1}}};\underline a)
\end{equation}
The point $[0:\dots:0:1]\in\mathbf P_A^{m+1}
=\mathbf P{\mathsf H}_{p^n}\otimes A$ 
is fixed and one easily checks 
that the action extends on the blow-up
$\tilde\mathbf P^{m+1}\subset
\mathbf P^{m+1}\times\mathbf P^m$; 
in particular, the inertia subgroup of the exceptional 
divisor for the action of 
$\mathbf Z/p^{n+1}\mathbf Z={\mathbb W}_{n+1}(\mathbf F_p)$
has order $p$. 
As $\overline{\mathbb W}_{n+1}=\tilde\mathbf P^{m+1}
\times_{\mathbf P^{m}}\overline{\mathbb W}_n$, 
we have the obtained the desired action;
moreover, the 
$\mathbf Z/p^{n+1}\mathbf Z$-inertia of 
the infinity section 
$\Sigma_{n+1}$ is of order $p$. 

\noindent
Let $\pi:P\to\overline{\mathbb W}_{n+1}$ 
be the cyclic cover of degree $p^n$
given by the fibred product
$$
\CD
P @>{\pi}>> \overline{\mathbb W}_{n+1} \\
@V{q}VV    @VV{\mathbf r}V \\
\overline{\mathbb W}_{n} @>{\mathbf\Psi_n}>> \overline{\mathbb W}_n
\endCD
$$
As $\mathbf\Psi_n^\ast$ is multiplication by $p$,
we have $P=\mathbf{P}\left(\mathcal{O}_{{\overline{\mathbb W}}_n}
\oplus\mathcal{O}_{{\overline{\mathbb W}}_n}(p^2)\right)$.
We have a sheaf of graded 
$\mathcal{O}_{{\overline{\mathbb W}}_n}$-algebras
$$
\mathcal S=\bigoplus_{d\geq0}
\mathcal{O}_{{\overline{\mathbb W}}_n}(pd)
$$
and, with the notations as in \cite{EGA2} \S2,
$\overline{\mathbb W}_{n+1}=\mathrm{Proj}\ \mathcal S[Y_n]$
and 
$P=\mathrm{Proj}\ \mathcal S^{(p)}[X_n]$.
Define a morphism of graded 
$\mathcal{O}_{{\overline{\mathbb W}}_n}$-algebras:
$$
\varphi^\#:\mathcal S^{(p)}[X_n]
\longrightarrow
\left(\mathcal S[Y_n]\right)^{(p)}
$$
which is the identity on $\mathcal S^{(p)}$ and sends $X_n$ to
$$
T^{p^{n+1}}S_n(\frac{Y_0^p}{T^p},\dots,
\frac{Y_n^p}{T^{p^{n+1}}};
-(\frac{Y_0}{T},\dots,\frac{Y_n}{T^{p^n}}))
=Y_n^p-Y_nT^{p^n(p-1)}
+T^{p^{n+1}}c_n(\frac{Y_0^p}{T^p},
\dots;-\frac{Y_0}{T},\dots).
$$
This defines a $\overline{\mathbb W}_n$-morphism 
$\varphi: G(\varphi^\#)\to P$
from an open subset of $\mathrm{Proj}\ (\mathcal
S[Y_n])^{(p)}=\overline{\mathbb W}_{n+1}$
(\textit{loc.cit}, prop. 3.1.8). 
The only homogeneous prime ideal containing
$\mathcal S^{(p)}$ and $Y_n$ is the irrelevant ideal,
hence $G(\varphi^\#)=\overline{\mathbb W}_{n+1}$.
As $\varphi$ is equivariant under the action
of $\mathbf Z/p\mathbf Z$ given by
equation (\ref{act}) (via the $n$-th 
iteration of the Verschiebung),
the composite morphism
$$
\CD
\mathbf\Psi_{n+1}:\overline{\mathbb W}_{n+1}
@>{\varphi}>> P @>{\pi}>>
\overline{\mathbb W}_{n+1}
\endCD
$$
extends the Artin-Schreier-Witt isogeny and is 
equivariant under the 
$\mathbf Z/p^{n+1}\mathbf Z$-action defined above
(it suffices to check it over the dense open subset 
${\mathbb W}_{n+1}$) and satisfies 
$$
\mathbf \Psi_n\circ\mathbf r
=\mathbf r\circ\mathbf \Psi_{n+1}.
$$
To prove the last statement, recall 
that, via $\mathbf r^\ast$, 
the cohomology group 
$A^\ast(\overline{\mathbb W}_{n+1})$
is a free
$A^\ast(\overline{\mathbb W}_{n})$-module
generated by the classes
$1$ and $\xi=[Z_{n+1}]$.
As $\mathbf\Psi$ commutes with $\mathbf r$,
it suffices to prove that
$\mathbf \Psi_{n+1}^\ast\xi=p\,\xi$.
Write:
\begin{equation}
\label{xi}
\mathbf\Psi_{n+1}^\ast\xi=
a\xi+\mathbf \beta \qquad
a\in\mathbf Z,\,\beta
\in{\mathbf r}^\ast A^1(\overline{\mathbb W}_{n}).
\end{equation}
Intersecting with a fibre of ${\mathbf r}$, we deduce $a=p$.
On the other hand, let 
$\eta={\mathbf r}^\ast[Z_n]$; 
since $c_t(\mathcal{O}_{{\overline{\mathbb W}}_n}
\oplus\mathcal{O}_{{\overline{\mathbb W}}_n}(p))
=c_t(\mathcal{O}_{{\overline{\mathbb W}}_n}(p))
=1+p\,[Z_n]t$, we have 
$$
\xi^2=p\,\xi\eta.
$$
As $\mathbf\Psi_{n+1}^\ast\eta=p\,\eta$, taking 
self intersection on both sides of (\ref{xi}) we get:
$$
p^2\xi^2+p^2\beta\eta=
p(p\xi+\beta)p\eta=
\mathbf\Psi_{n+1}^\ast(p\xi\eta)
=\mathbf\Psi_{n+1}^\ast(\xi^2)
=p^2\xi^2+2p\xi\beta+\beta^2
$$
whence $(2p\beta)\,\xi+(\beta^2-p^2\beta\eta)=0$.
Since $1,\xi$ form a 
$A^\ast(\overline{\mathbb W}_{n})$-basis,
this implies $\beta=0$.

\begin{cor}
With the notation as in Corollary~\ref{bord}, 
the inertia subgroup at the 
$i$-th component $B_{n,i}$ of the boundary 
divisor $B_n$ on $\overline{\mathbb W}_{n}$
is cyclic of order $p^{n-i}$.
\end{cor}

Let now $f_n:C_n\to D$ be a cyclic cover as 
in Theorem 1 and 
$\underline u: D-\{x_0\}\to {\mathbb W}_{n}$ the 
corresponding Witt representative; 
clearly, $\underline u$ extends to
$$
\underline u: D\longrightarrow \overline{\mathbb W}_{n}.
$$

\begin{lemma} 
With the notation as above, $\deg\underline u^\ast
\mathcal{O}_{{\overline{\mathbb W}}_n}(1)=
\max_{0\leq i\leq n-1}\{p^{n-i-1}\nu_i\}$.
\end{lemma}
Proof:
It suffices to consider the composite map 
(cf. lemma~\ref{H=Gamma})
$$
\CD
D@>{\underline u}>> {\overline{\mathbb W}}_n
@>{\phi_{|Z_n|}}>> \mathbf P{\mathsf H}_{p^{n-1}}.
\endCD
$$
In a neighborhood of $x_0$, it is defined by
$[s^{M}:s^{M}u_0^{p^{n-1}}:\dots,s^{M}u_{n-1}]$, 
where $s$ is a local parameter at $x_0$ and 
$$
\begin{array}{rl}
M & =\max\,\{i_0\nu_0+\dots+i_{n-1}\nu_{n-1}|\,
0\leq i_h\leq p^{n-1-h}\,
\mathrm{and}\,
\sum_{h=0}^{n-1} p^hi_h=p^{n-1}\}\\
& = \max_{0\leq i\leq n-1}\{p^{n-i-1}\nu_i\}.
\end{array}
$$

\begin{prop}
\label{bound}
Let $f_n:C_n\to D$ be as in Theorem 1 and let 
$M(f_n)=\deg\underline{u}^\ast\mathcal{O}_{{\overline{\mathbb W}}_n}(1)$.
If $m(f_n)+1$ is the degree of the conductor, then
$$
m(f_n)\leq M(f_n).
$$
\end{prop}

\noindent
Proof: We show that $\mathfrak{u}=(1+M(f_n)).\infty$ is a modulus, in the
sense of~\cite{gacc}, chap. III, for the rational
map $\underline{u}:D-\{\infty\}\to\Bbb{W}_n$, with the smallest possible degree.
By the universal property of the generalized jacobians, $\underline{u}$ factors
through $J_{\mathfrak{u}}$ and, by the minimality of the conductor,
$\mathfrak{m}\leq\mathfrak{u}$.

We have to show that 
$(\underline{u},\alpha)_\infty=\underline{0}\in{\mathbb W}_n(k)$
for any rational function $\alpha$ on
$D$ such that $\alpha\equiv1 \bmod \mathfrak{u}$, where
$(\underline{u},-)_\infty$ denotes the local symbol associated to 
$\underline{u}$: it is the "Residuenvektor" defined in~\cite{witt}, \S 2, 
whose $j$-th phantom component is 
$$
\Phi_j\left((\underline{u},\alpha)_\infty\right)
=\mathrm{Res}\,\left(\Phi_j(u_0,\dots,u_j)\frac{d\alpha}\alpha\right)
$$
Let $s$ be a local parameter on $D$ at $\infty$; if $A$ is any complete 
discrete valuation ring of characteristic zero with residue field $k$, 
we can lift $\alpha\in k[[s]]$ and $u_0,\dots,u_{n-1}\in k((s))$
to formal power series $\tilde{\alpha}\in A[[s]]$ and
$\tilde{u}_0,\dots,\tilde{u}_{n-1}\in A((s))$, with the
$\tilde{u}_i$ having a pole at $s=0$ of the same order as $u_i$:
by definition, $(\underline{u},\alpha)$ is the reduction  modulo $pA$ of
$(\underline{\tilde u},\tilde\alpha)$. 

\noindent
Recall that $u_i$ has a pole of order $\nu_i$ and that 
$M=\max\,\{p^{n-i-1}\nu_i:\,0\leq i\leq n-1\}$. 
It is then clear that $M$ is the smallest integer $\rho$ for which 
$$
\begin{array}{rl}
\Phi_j\left((\underline{\tilde u},\tilde \alpha)\right)
& = 
\mathrm{Res}\,\left(\Phi_j(\tilde{u}_0,\dots,\tilde{u}_j)
\frac{d\tilde\alpha}{\tilde\alpha}\right) \\ 
& =
\mathrm{Res}\,\left(\left(\tilde{u}_0^{p^j}+\dots+p^j\tilde{u}_j\right)
\frac{d\tilde\alpha}{\tilde\alpha}\right)
\end{array}
$$
vanishes for any $0\leq j\leq n-1$ and for any $\tilde\alpha$
such that $1-\tilde\alpha$ has a zero of order $\rho+1$ at $s=0$.

\bigskip

To complete the proof of the Theorem, 
we need to show that the bound
provided by Proposition~\ref{bound} is attained.
When $k$ is finite, this was done by H.L. Schmid, 
\cite{schmid}, \S3:
if $\chi:G_{k(D)}^{ab}\to \mathbf Z/p^n\mathbf Z$ 
is the character corresponding to the cover $f_n$,
then the residue vector is related to the reciprocity 
map by the formula
$$
\chi(\mathrm{rec}(\alpha))=
\mathrm{Tr}\left((\underline u,\alpha)_\infty\right)
\qquad \forall\alpha\in k(D)^\times,
$$
where $\mathrm{Tr}:{\mathbb W}_n(k)\to{\mathbb W}_n(\mathbf F_p)$
is the trace map.
The computation of the conductor follows at once
from its classical description in terms of norms of 
units (e.g. \cite{CL}, chap. XV, \S2).

For an arbitrary field $k$, this classical approach no 
longer holds, but it is still possible to prove the desired
equality using higher class field theory: see \cite{bry},
\S2.

In the appendix we give an elementary proof of this formula;
the reader who is willing to venture through it will 
certainly appreciate the difference between the clean 
and conceptual approach via class field theory and the 
following tedious but elementary computations.

\bigskip
\bigskip
\bigskip

\appendix

\noindent
{\bf Appendix: Elementary computation of the 
conductor of a cyclic extension.}

\bigskip
\bigskip

Our goal is to show that the inequality in
Proposition~\ref{bound} is indeed an equality. 
Again, our proof is by induction on $n$, the case $n=1$ being
classical (e.g.~\cite{clac}, 4.4). We shall therefore assume that the
Theorem is true for cyclic extensions of degree $p^l$, for $l\leq n$ and
prove that it is true for $l=n+1$ by showing that, with the notations of
Proposition~\ref{bound}, $m(f_{n+1})=M(f_{n+1})$. 

To this end, we need to compute the conductor of $C_{n+1}/D$ or, which
amounts to the same, that of $C_{n+1}/C_n$.
The latter is a cyclic extension of degree $p$ and hence it is described
by an  Artin-Schreier equation $y^p-y=z$, for a suitable rational function
$z$ on $C_n$ which is not, however, in standard form i.e. $z$ has a pole
of order divisible by $p$, and we cannot read off the conductor
immediately.

\bigskip

Some work is needed in order to bring the Artin-Schreier equation
describing the cover $C_{n+1}/C_n$ in standard form. 
The key to all our computations is the following, elementary

\begin{lemma} 
\label{adjust}
(Adjustment lemma).
Let $f:Y\to X$ be a finite separable cover of curves of degree
$p^n$, totally ramified at $y_0$. Let $t$ (resp. $s$) be local parameters
at $y_0$ (resp. $x_0=f(y_0)$). Then, in a neighborhood of $y_0$
$$
s=t^{p^n}\left(\alpha_1^p+t^{\mu(f)}\alpha_0\right)
$$
for some $\alpha_0,\alpha_1\in\mathcal{O}_{Y,y_0}^\times$, with
$$
\mu(f)=\deg\mathfrak{D}_{Y/X,y_0}-p^n+1
$$
\end{lemma}

\noindent
Proof: Let $s=\alpha t^{p^n}$ for some unit
$\alpha\in\mathcal{O}_{Y,y_0}^\times$ 
and expand $\alpha=\sum_{r\geq0}a_rt^r$ with
$a_r\in k$ and $a_0\in k^\times$. 
If $\mu$ is the smallest positive integer prime to $p$ such that
$a_\mu\neq0$, then 
$$
\deg\mathfrak{D}_{Y/X,y_0}=v_{y_0}\left(\frac{ds}{dt}\right)
=p^n+v_{y_0}\left(\frac{d\alpha}{dt}\right)=p^n+\mu-1
$$
We can therefore collect:
$\alpha=\sum_{0\leq r<\frac\mu p}a_{pr}t^{pr}+t^\mu\sum_{r\geq\mu}a_rt^r
=\alpha_1^p+t^\mu\alpha_0$.

\begin{cor}
With the notations above, let $g$ be a rational function on $X$. 
Then $f^\ast(g)=g_1^p+g_0$ for suitable rational
functions $g_i$ on $Y$ with $v_{y_0}(g_0)=p^nv_{x_0}(g)+\mu(f)$.
\end{cor}

\noindent
Proof: Expanding $g$ as a power series in $s$, it is sufficient to prove
it when $g(s)=as^r$ is a monomial of degree prime to $p$:
$$
as^r=at^{p^nr}\left(\alpha_1^p+t^\mu\alpha_0\right)^r
=at^{p^nr}\left(\alpha_1^{pr}+t^\mu\alpha_0^\prime\right)
=(a^\frac1pt^{p^{n-1}}\alpha_1^r)^p+at^{p^nr+\mu}\alpha_0^\prime
$$

\bigskip

When the cover is Galois, it is easy to relate the
number $\mu(f)$ to the conductor of the extension and how it varies in
towers:

\begin{prop}
\label{condor}
With the notation as in Lemma~\ref{adjust}, suppose moreover that $f$ is
Galois of group $G$ and conductor $m(f)+1$. 
\begin{enumerate}
\item
If $e(f)$ is the largest
integer such that $G_e\neq1$, then
$$
\mu(f)=p^nm(f)-e(f)
$$
\item
Let $g:Z\to Y\to X$ be a second separable Galois cover, 
totally ramified above $x_0$ whose group
$\tilde G$ is an extension of $G$ by $H=\mathbf{Z}/p\mathbf{Z}$. 
Let $e(g)$ be the largest integer such that $\tilde G_e\neq1$ 
and suppose that $\tilde G_{e(g)}\cap H\neq1$.
If $m(g)+1$ denotes the conductor of $g$, then
$$
e(g)=p^nm(g)-\mu(f)
$$
\end{enumerate}
\end{prop}

\noindent
Proof: (1). Denote, as usual, by $g_i=|G_i|$; 
as $f$ is totally ramified,
$G_0=G$. By\cite{CL}, chap. IV, proposition 12,
$m=\frac1{g_0}\sum_{i=0}^eg_i-1$, hence
$$
\begin{array}{rl}
\mu(f) & =\deg\mathfrak{D}_{Y/X,y_0}-p^n+1 \\
& = \sum_{i=0}^e(g_i-1)-p^n+1 \\
& = \sum_{i=0}^eg_i-e-p^n \\
& = p^n\left(\frac1{p^n}\sum_{i=0}^eg_i-1\right)-e.
\end{array}
$$
(2). Denote by $\varphi$ the Herbrand function of the
cover $g$ and by $\tilde g_i=|\tilde G_i|$. 
The hypotheses imply that $\tilde g_i=pg_i$ 
for $0\leq i\leq e(f)$, hence
$$
\begin{array}{rl}
m(g)=\varphi(e(g)) &
=\frac1{p^{n+1}}\left(\tilde g_1+\dots+\tilde g_{e(f)}+
p\,[e(g)-e(f)]\right) \\  
& = \frac1p\varphi(e(f))+\frac1{p^n}[(e(g)-e(f)] \\
& = m(f)+\frac1{p^n}[e(g)-e(f)] .
\end{array}
$$

\begin{cor}
\label{tele}
Let $f_n:C_n\to D$ be a separable cyclic cover of degree $p^n$, totally
ramified above $\infty$. 
For any $0\leq i\leq n$, let $f_i:C_i\to D$ be the quotient cover of
degree $p^i$ and by $g_i:C_n\to C_i$ the subcover of degree $p^{n-i}$.
Put $m_i=m(f_i)$ (resp. $e_i=e(f_i)$, resp. $\mu_i=\mu(f_i)$). Then:
$$
e_n=\sum_{i=1}^np^{i-1}(m_i-m_{i-1}); \qquad
\mu_n=\sum_{i=1}^n(p^i-p^{i-1})m_i; \qquad
\mu(g_i)=\mu_n-p^{n-i}\mu_i.
$$
\end{cor}

\noindent
Proof: From Proposition~\ref{condor}.2, we deduce the inductive relations 
$$
e_n-e_{n-1}=p^{n-1}(m_n-m_{n-1});\qquad
\mu_n-\mu_{n-1}=(p^n-p^{n-1})m_n$$
from which the first two formulas in the claim follow at once. 
As for the third, it is a consequence of the transitivity of the different
of a composite extension:
$$
\mu(g_i)=\deg\mathfrak{D}_{C_n/C_i}-p^{n-i}+1
=\deg\mathfrak{D}_{C_n/D}-p^n+1
-p^{n-i}\left(\deg\mathfrak{D}_{C_i/D}-p^i+1\right).
$$

\bigskip

We are now ready to resume the discussion interrupted before
Lemma~\ref{adjust}. 
We are given a separable cyclic cover $f_{n+1}:C_{n+1}\to D$ of 
degree $p^{n+1}$, as in the theorem, and we want to compute 
its conductor, to show that the inequality provided by 
Proposition~\ref{bound} is indeed an equality.

Equivalently, we can compute the conductor of the  degree $p$
subextension $C_{n+1}/C_n$. 
With the notations as in Proposition~\ref{condor}, 
let $\tilde G=\mathbf{Z}/p^{n+1}\mathbf{Z}$ be the group of 
$C_{n+1}/D$ and $H=\mathbf{Z}/p\mathbf{Z}$ that of 
$C_{n+1}/C_n$:  
by~\cite{CL}, IV, prop. 2, $H_i=H\cap\tilde G_i$
and it is clear that $H=H_{e_{n+1}}=\tilde G_{e_{n+1}}$ is the last
nontrivial subgroup for both extensions.

\bigskip

We begin by making the Artin-Schreier equation
defining the cover $C_{n+1}/C_n$ more explicit.
Away from $\infty$, the cover $C_{n+1}/D$ is described by the Witt 
vector equation
\begin{equation}
\label{WAS}
F(\underline{Y})-\underline{Y}=(u_0,\dots,u_n)
\end{equation}
where the left hand side should be understood 
as a difference of two
vectors, while the equation defining $C_n/D$ 
is the image of (\ref{WAS})
under the restriction homomorphism 
$\mathbf r:\Bbb{W}_{n+1}\to\Bbb{W}_n$.
Therefore, the $n$-th component of (\ref{WAS}) is
$$
Y_n^p-Y_n+c_n\left(F(\underline{Y}),-\underline{Y}\right)=u_n
$$
where the polynomial $c_n(\underline{X},\underline{Y})$ 
is the one defined by equation (\ref{S_n}). 
If $\underline{y}=(y_0,\dots,y_{n-1})$ denotes
any solution of the Artin-Schreier-Witt 
equation defining $C_n/D$, the
cover
$C_{n+1}/C_n$ is described by the Artin-Schreier equation
\begin{equation}
\label{AS}
Y_n^p-Y_n=u_n-c_n\left(\underline{y}^p,-\underline{y}\right).
\end{equation}

As remarked, this equation is not in standard form: the datum on the
right hand side has, in general, poles of order divisible by $p$. 
We can get rid of them by adding terms of the type 
$\wp(h)=h^p-h$ for suitable rational functions $h$ on $C_n$; the result
will be a rational function with a pole at $\infty$ of order 
$e_{n+1}$. 
By Proposition~\ref{condor}.2, $e_{n+1}=p^nm_{n+1}-\mu_n$.

To complete the proof of the theorem, we should then prove the 
following:

\begin{claim}
The right hand side of equation (\ref{AS}) is congruent
$\bmod\,\wp\left(k(C_n)\right)$ to a function with 
a pole at $\infty$ of order $p^nM-\mu_n$, where 
$M=\max\,\{p^{n-i}\nu_i:\,0\leq i\leq n\}$.
\end{claim}

It is easy to handle the term $u_n$ in (\ref{AS})
with the adjustment lemma:
$u_n$ is a rational function on $D$ 
with a pole of order $\nu_n$ at
$\infty$, hence we can find rational functions
$u_n^\prime,\,u_n^{\prime\prime}$ on $C_n$ such that:
\begin{equation}
\label{new}
u_n={u_n^\prime}^p+u_n^{\prime\prime}
\qquad \mathrm{with}\
v_n(u_n^{\prime\prime})=-p^{n}\nu_n+\mu_n
\end{equation}
We can then rewrite 
$u_n=u_n^{\prime\prime}+u_n^\prime+\wp(u_n^\prime)$.
When $M=\nu_n$, the function $u_n^{\prime\prime}$ 
has the predicted valuation and in any event 
the contribution of $u_n^\prime$ is negligible, as
$v_n(u_n^\prime)=-p^{n-1}\nu_n\geq-p^{n-1}M$, hence
\begin{equation}
\label{rough}
\begin{array}{rl}
-p^nM+\mu_n+p^{n-1}\nu_n &\leq
-(p^n-p^{n-1})M+\sum_{i=1}^n(p^i-p^{i-1})m_i\\
& < -m_n(p^{n+1}-2p^n+1) <0
\end{array}
\end{equation}
as $m_i<m_n$ for $i<n$ and $M\geq pm_n$ (by our inductive hypothesis, 
$m_{j+1}=\max\,\{p^{j-i}\nu_i:\,i\leq j\}$ for $j\leq n-1$).

\bigskip

We must then analyze the contribution of the second term on the right
hand side of (\ref{AS}); according to the claim, it should be dominant 
when $\nu_n<M$. Indeed, we will prove the following:

\begin{claim}
\label{old}
The function
$c_n(\underline{y}^p,-\underline{y})$ is congruent
$\bmod\,\wp\left(k(C_n)\right)$ to a function with 
a pole at $\infty$ of order $p^{n+1}m_n-\mu_n$.
\end{claim}

By the induction hypothesis, $M=\max\,\{\nu_n,pm_n\}$, and these two
numbers cannot be equal; hence, the latter claim implies the former.

The following lemma gathers the information we shall need about the
polynomials $c_i$.

\begin{lemma}
\label{c_n}
Fix integers $n$ and $i\leq n-1$. Then
$$
c_n(\underline{X},\underline{Y})=-X_i^{p^{n-i}-1}
\left(Y_i+c_i(\underline{X},\underline{Y})\right)
+R_{i,n}
$$
where $R_{i,n}$ is an isobaric polynomial of weight $p^n$ in which the
variable
$X_i$ appears with degree strictly smaller than $p^{n-i}-1$.
\end{lemma}

\noindent
Proof: It is well known that $c_0=0$ and
$c_1=-\sum_{r=1}^{p-1}\frac1p{p\choose r}X_0^rY_0^{p-r}$. 
We proceed by induction. Using the recursive relation 
$\Phi_n(X_0,\dots,X_n)=\Phi_{n-1}(X_0^p,\dots,X_{n-1}^p)+p^nX_n$,
we see that
$$
c_n=S_n-X_n-Y_n
=\frac1{p^n}\left[\sum_{h=0}^{n-1}p^h\left(X_h^{p^{n-h}}+Y_h^{p^{n-h}}
-S_h^{p^{n-h}}\right)\right].
$$
Fix $0\leq h\leq n-1$; applying the binomial formula twice:
\begin{equation}
\begin{array} {l}
X_h^{p^{n-h}}+Y_h^{p^{n-h}}
-\left(X_h+Y_h+c_h\right)^{p^{n-h}}=\\
\mbox{}\\
-c_h^{p^{n-h}} 
-\sum_{s=1}^{p^{n-h}-1}{p^{n-h}\choose s}
X_h^sY_h^{p^{n-h}-s}-
\sum_{r=1}^{p^{n-h}-1}\sum_{s=0}^{r}
{p^{n-h}\choose r}{r\choose s}X_h^sY_h^{r-s}c_h^{p^{n-h}-r}.
\end{array}
\label{h-term}
\end{equation}
The variable $X_i$ appears in the expression above only if $h\geq i$;
suppose $h>i$, then by the induction hypothesis we can write
$$
c_h=-X_i^{p^{h-i}-1}
\left(Y_i+c_i\right)
+R_{i,h}
$$
hence the highest power of $X_i$ in (\ref{h-term}) appears in the term
$c_h^{p^{n-h}}$ and its contribution in $c_n$ is of degree at most 
$(p^{h-i}-1)p^{n-h}$.

\noindent
For $h=i$, the variable $X_i$ does not appear in $c_i$ hence the
term containing the highest power of $X_i$ in (\ref{h-term}) is
$$
-{p^{n-i}\choose1}X_i^{p^{n-i}-1}Y_i-{p^{n-i}\choose
p^{n-i}-1}{p^{n-i}-1\choose p^{n-i}-1}X_i^{p^{n-i}-1}c_i
=-p^{n-i}X_i^{p^{n-i}-1}
\left(Y_i+c_i\right).
$$

\bigskip

Our task is now to compute the valuation 
at $\infty$ of a rational function on $C_n$, 
congruent to $c_n(\underline{y}^p,-\underline{y})
\bmod\,\wp(k(C_n))$ and having a pole of order prime to $p$,
where $\underline{y}=(y_0,\dots,y_{n-1})$ is a vector 
of solutions for the Artin-Schreier-Witt equation defining $C_n/D$.

Lest the reader should get lost or depressed (or both) by the 
following computations, let us state immediately that we are 
going to show that, for any $i\leq n-1$,  only the monomials 
$$
y_i^{p(p^{n-i}-1)+1}
$$
are relevant. More precisely, if $j\leq n-1$ is the 
(necessarily unique) integer for which
$$
m_n=p^{n-j-1}\nu_j=\max\,\{p^{n-i-1}\nu_i:\,0\leq i\leq n-1\}
$$
then $y_j^{p(p^{n-j}-1)+1}$ is congruent
$\bmod\,\wp\left(k(C_n)\right)$ to a function with 
a pole at $\infty$ of order $p^{2n-j}\nu_j-\mu_n$ 
and that 
$c_n(\underline{y}^p,-\underline{y})-y_j^{p(p^{n-j}-1)+1}$ is
congruent  to a function with a
pole of order strictly smaller. 

\bigskip

The approach is in two steps: first, we obtain an 
estimate of the valuation of 
$c_n(\underline{y}^p,-\underline{y})$
and then we apply the adjustment lemma to the 
functions $y_i$ to get rid of the terms in 
$\wp(k(C_n))$.
In both steps, the problem to sort out the relevant monomials is reduced
to a linear optimization.

The proof is only a tedious exercise in linear programming, and the
reader would loose very little, should he indulge in the temptation 
to skip the rest of this section. 

\bigskip

When the variables $X_i$ and $Y_i$ are 
given the weight $p^i$,
the polynomial $c_n(\underline{X},\underline{Y})$ 
is isobaric of weight $p^n$; evaluating at
$\underline{X}=\underline{y}^p$ and 
$\underline{Y}=-\underline{y}$, our task is to compute 
the valuation of a sum of monomials $\prod_{i=0}^{n-1}y_i^{pa_i+b_i}$.

Put $v_n(y_i)=w_i$; regarding $a_i$ and $b_i$ ad variables,
the question can be rephrased as an optimization problem for a linear
function under a linear constraint:
\begin{equation}
\label{LP1}
\mathrm{Minimize:} \sum_{i=0}^{n-1}w_i(pa_i+b_i)\qquad
\mathrm{under:}\,
\sum_{i=0}^{n-1}p^i(a_i+b_i)=p^n
\end{equation}
with $a_i\in[0,p^{n-i}-1]$ and $b_i\in[0,p^{n-i}]$.
The range of the $a_i$'s is restricted by Lemma~\ref{c_n}; 
we do not apply it to the $b_i$'s because of the irregular 
behavior in characteristic $2$ of the polynomials 
$I_i(\underline Y)$ expressing the opposite of a Witt vector.

Actually, it is more convenient to
rephrase the problem slightly:  
putting $\alpha_i=pa_i+b_i$ and taking
$\alpha_i,b_i$  as variables, it becomes:
\begin{equation}
\label{LP2}
\mathrm{Minimize:}\ 
\sum_{i=0}^{n-1}w_i\alpha_i\qquad
\mathrm{under:}\,
\sum_{i=0}^{n-1}p^i\alpha_i+\sum_{i=0}^{n-1}(p^{i+1}-p^i)b_i=p^{n+1}
\end{equation}
with $\alpha_i\in[0,p^{n-i+1}-p+1]$ and $b_i\in[0,p^{n-i}]$.

Extremal values of a linear
function on  a convex polytope are 
attained at the vertices;
as $b_i=p^{n-i}$ forces $a_i=0$, hence $\alpha_i=p^{n-i}$, 
best admissible solutions are for 
$\alpha_i=p^{n-i+1}-p+1$, $b_i=1$ and and $\alpha_l=b_l=0$ 
for $l\neq i$. We have thus:
$$
v_n(c_n(\underline{y}^p,-\underline{y}))
\geq\min\left\{(p^{n-i+1}-p+1)w_i\,|\quad0\leq
i\leq n-1\right\}
$$
and this in turn implies our first estimate:
\begin{lemma}
\label{sort}
\begin{enumerate}
\item
$v_n(c_n(\underline{y}^p,-\underline{y}))\geq-(p^{n+1}-p+1)m_n$.
\item
$v_{n+1}(y_n)\geq-p^nm_{n+1}$, with an equality if and only if 
$\nu_n=m_{n+1}$.
\end{enumerate}
\end{lemma}

\noindent
Proof: Induction on $n$: the case $n=0$ being clear,
assume (1) and (2) hold for 
$c_i$ and $y_i$ for all $i\leq n-1$.
Then, since $m_n\geq p^{n-i-1}m_{i+1}$, we have
$$
\begin{array}{rl}
v_n(y_i^{p(p^{n-i}-1)+1}) & = 
p^{n-i-1}v_{i+1}(y_i^{p(p^{n-i}-1)+1}) \\
& \geq
-p^{n-i-1}(p^{n+1}-p^{i+1}+p^i)m_{i+1} \\
& \geq 
-(p^{n+1}-p^{i+1}+p^i)m_n
\end{array}
$$
Finally, as $y_n$ is a solution of the equation 
$y_n^p-y_n=u_n-c_n$ and  
clearly $v_n(u_n)=-p^n\nu_n\geq-p^nm_{n+1}$, 
with equality holding if and only if 
$\nu_i=m_{i+1}$, (2) follows at once.

\bigskip

To compute the valuation of
$c_n(\underline{y}^p,-\underline{y})\bmod\wp(k(C_n))$, 
we have to apply the adjustment lemma twice, 
first on $C_i$  and then on $C_n$.

Fix $i\leq n-1$ and denote by $v_i$ the discrete valuation of
$\mathcal{O}_{C_i,\infty}$. 
By definition, $y_i$ is the solution (on
$C_{i+1}$)  of the Artin-Schreier equation
$y_i^p-y_i=u_i-c_i(\underline{y}^p,-\underline{y})$; 
applying the adjustment lemma to $u_i$ and 
the induction hypothesis to $c_i$, we can find rational functions $g_i$
and $\tilde u_i$ on $C_i$ such that 
$$
u_i-c_i(\underline{y}^p,-\underline{y})=\tilde u_i+g_i^p-g_i
$$
with $v_i(\tilde u_i)=p^{i}m_{i+1}-\mu_{i}$ and, by the lemma above,
$v_i(g_i^p)\geq-p^nm_{n+1}$, equality holding iff $\nu_i=m_{i+1}$.

On $C_{i+1}$ we have thus $y_i=\tilde y_i+g_i$ with
$$
v_{i+1}(\tilde y_i)=-p^im_{i+1}+\mu_i;
\qquad
v_{i+1}(y_i)=v_{i+1}(g_i)\geq p^im_{i+1}
$$
with an equality iff $\nu_i=m_{i+1}$.

We now apply the adjustment lemma on $C_n$ both to
$\tilde y_i$ and $g_i$: we can write
$$
\tilde y_i=\tilde\gamma_i^p+\tilde\eta_i \qquad
g_i={\gamma_i^\prime}^p+\eta_i^\prime
$$
with
$v_n(\tilde\gamma_i)=p^{-1}v_n(y_i)=-p^{n-2}m_{i+1}+p^{n-i-1}\mu_{i+1}$
and, applying Corollary~\ref{tele},
$$
\begin{array}{rl}
v_n(\tilde\eta_i) & =v_n(\tilde y_i)+\mu(C_n/C_{i+1}) \\
& = -p^{n-i-1}(p^im_{i+1}-\mu_i)+ \mu_n-p^{n-i-1}\mu_{i+1} \\
& = -p^{n-1}m_{i+1}-p^{n-i-1}(p^{i+1}-p^i)m_{i+1}+\mu_n \\
& = -p^nm_{i+1}+\mu_n.
\end{array}
$$

As $g_i$ comes from $C_i$, we have
$v_n(\gamma_i^\prime)=p^{-1}v_n(g_i)\geq-p^{n-2}m_{i+1}$, with equality
iff $\nu_i=m_{i+1}$ and 
$$
\begin{array}{rl}
v_n(\eta_i^\prime) & =v_n(g_i)+\mu(C_n/C_i) \\
& = p^{n-i}v_i(g_i)+\mu_n-p^{n-i}\mu_i \\
& \geq -p^{n-1}m_{i+1}+\mu_n-p^{n-i}\mu_i
\end{array}
$$
Trivially, $v_n(\gamma_i^\prime)\lneq v_n(\tilde\gamma_i)$; the same rough
estimates as in equation (\ref{rough}) yield 
$v_n(\tilde\eta_i)\lneq v_n(\eta_i^\prime)$. Putting 
$\gamma_i=\tilde\gamma_i+\gamma_i^\prime$ and 
$\eta_i=\tilde\eta_i+\eta_i^\prime$, we can write
$$
y_i=\gamma_i^p+\eta_i \qquad \textrm{with:} \quad
v_n(\eta_i)=-p^nm_{i+1}+\mu_n;
\ v_n(\gamma_i)\geq-p^{n-2}m_{i+1}
$$
again the last inequality being an equality iff $\nu_i=m_{i+1}$.

Substituting in $\prod_{i=0}^{n-1}y_i^{\alpha_i}$ we get:
$$
\prod_{i=0}^{n-1}(\gamma_i^p+\eta_i)^{\alpha_i}
=\prod_{i=0}^{n-1}\gamma_i^{p\alpha_i}+ 
\sum_{i=0}^{n-1}\left(\alpha_i\gamma_i^{p\alpha_i-1}\eta_i
\prod_{h\neq i}\gamma_h^{p\alpha_h}\right)+
\dots
$$
The first term is a $p$-power and, by Lemma~\ref{sort}, 
its valuation is bounded from below by $-(p^{n+1}-p+1)m_n$; 
the same estimates as those of formula (\ref{rough})
show that its $p$-th root has a pole of order strictly 
smaller than $p^{n+1}m_n-\mu_n$ and we can get rid of it 
$\bmod\wp(k(C_n))$.

\bigskip

To prove the claim, we should therefore compute, 
for all $i,h\leq n-1$, the valuation of the monomials
$\gamma_i^{p\alpha_i-1}\eta_i
\prod_{h\neq i}\gamma_h^{p\alpha_h}$.

Our optimization problem (\ref{LP2}) becomes then:
$$
\mathrm{Minimize:}\ v_n(\eta_i)+w_i\alpha_i+
\sum_{h\neq i}w_h\alpha_h\qquad
\mathrm{under:}\,
\sum_{l=0}^{n-1}p^l\alpha_l+\sum_{l=0}^{n-1}(p^{l+1}-p^l)b_l=p^{n+1}
$$
with $\alpha_l\in[0,p^{n-l+1}-p+1]$ and $b_l\in[0,p^{n-l}]$.

\bigskip

\noindent
\textbf{Proof of Claim \ref{old}.}
All we have to do is check the extremal values of the function above,
for all $i$.
For $\alpha_i=p^{n-i+1}-p+1$ and $b_i=1$ we get
$$
v_n(\gamma_i^{p^{n-i+2}-p^2}\eta_i)\geq
-(p^{n-i+2}-p^2)p^{n-2}m_{i+1}-p^nm_{i+1}+\mu_n=-p^{2n-i}m_{i+1}+\mu_n
$$
with an equality iff $\nu_i=m_{i+1}$. 
Hence, if $j\leq n-1$ is the 
(necessarily unique) integer for which
$$
m_n=p^{n-j-1}\nu_j=\max\,\{p^{n-i-1}\nu_i:\,0\leq i\leq n-1\}
$$
then $y_j^{p(p^{n-j}-1)+1}$ is congruent
$\bmod\,\wp\left(k(C_n)\right)$ to a function with 
a pole at $\infty$ of order $p^{2n-j}\nu_j-\mu_n$. 
If $i\neq j$ then $p^{n-i-1}\nu_i\lneq m_n$
and the contribution of $y_i$ is strictly smaller. 

\bigskip

\noindent
Finally, consider
$c_n(\underline{y}^p,-\underline{y})-y_j^{p(p^{n-j}-1)+1}$: 
repeating the argument above,
to show that it is congruent $\bmod\wp(k(C_n))$ 
to a function with a pole of order
strictly smaller than $p^{2n-j}\nu_j-\mu_n$, 
it suffices to show that the term containing the 
highest power of $y_j$, namely
$$
y_j^{p(p^{n-j}-1)}c_j(\underline{y}^p,-\underline{y})
$$
(cf. Lemma \ref{c_n}) satisfies this property.
As above, after a double application of the adjustment 
lemma on $C_n$, we can write 
$$
c_j(\underline{y}^p,-\underline{y})=\theta_j^p+\xi_j
$$
with 
$v_n(\theta_j^p)=p^{n-j}v_j(c_j)\geq-p^{n-j}(p^{j+1}-p+1)m_j$ 
(cf. Lemma \ref{sort}) and 
$$
v_n(\xi_j)=-p^{n-j}(p^{j+1}m_j-\mu_j)+\mu_n-p^{n-j}\mu_j
=-p^{n+1}m_j+\mu_n.
$$
Therefore, $\bmod\wp(k(C_n))$ we can get rid of
$y_j^{p(p^{n-j}-1)}\theta_j^p$ and, by Lemma \ref{sort}.2, 
$$
\begin{array}{rl}
v_n(y_j^{p(p^{n-j}-1)}\xi_j) & =
p(p^{n-j}-1)p^{n-j-1}v_{j+1}(y_j)
-p^{n+1}m_j+\mu_n\\
& = -(p^{n+1}-p^{j+1})p^{n-j-1}\nu_j
-p^{n+1}m_j+\mu_n\\
& = -p^{2n-j}\nu_j+p^j(\nu_j-pm_j)+\mu_n
\end{array}
$$
and the assumption on $j$ implies that $\nu_j-pm_j$ is strictly
positive.

\end{document}